\documentclass{article}

\usepackage[all]{xy}
\usepackage[usenames,dvipsnames]{pstricks}
\usepackage{epsfig}
\usepackage{pst-grad} 
\usepackage{pst-plot} 
\usepackage{ebezier}
\usepackage{latexsym}\usepackage{amsfonts}\usepackage{amsthm}
\usepackage{color}
\usepackage{amssymb, amsmath}
\usepackage{amsthm}
\usepackage{newlfont}\usepackage{graphicx}\usepackage{doc}
\usepackage{mathrsfs}

\usepackage[
        colorlinks, citecolor=red,
        backref,
        pdfauthor={MDJ,JK},
        pdftitle={S1}
]{hyperref}

\newtheorem{thrm}{Theorem}
\newtheorem{lem}[thrm]{Lemma}

\newtheorem{defn}[thrm]{Definition}

\begin{document}

\title{Shadowing property on hyperspace of continua induced by Morse gradient system}

\author{
Jelena Kati\'c\thanks{Corresponding author: Jelena Kati\'c.}\\
Matemati\v cki fakultet\\
Studentski trg 16\\
11000 Beograd\\
Serbia\\
jelena.katic@matf.bg.ac.rs
\and
Darko Milinkovi\'c
\thanks{The work of both authors is partially supported by the Ministry of Education, Science and Technological Developments of Republic of Serbia: grant number  451-03-47/2023-01/ 200104 with Faculty of Mathematics.} \\
Matemati\v cki fakultet\\
Studentski trg 16\\
11000 Beograd\\
Serbia\\
darko.milinkovic@matf.bg.ac.rs
}

\maketitle

\begin{abstract} It is known that Morse-Smale diffeomorphisms have the shadowing property; however, the question of whether $C(f)$ also has the shadowing property when $f$ is Morse-Smale remains open and has been resolved only in a few specific cases~\cite{AB}. We prove that if $f:M\to M$ is a time-one-map of Morse gradient flow, the induced map $C(f):C(M)\to C(M)$ on the hyperspace of subcontinua does not have the shadowing property. 
\end{abstract}

\medskip


\section{Introduction}

A dynamical system is said to have the {\it shadowing property} (also known as the {\it pseudo-orbit tracing property}) if, informally speaking, every approximate orbit with small errors (i.e., a pseudo-orbit) can be closely followed by a true orbit. This concept was originally studied by  Anosov~\cite{A}, Bowen~\cite{Bo} and Sinai~\cite{Si}. If a dynamical system undergoes a small perturbation, the orbits of the perturbed system become pseudo-orbits of the original one. Therefore, shadowing is closely related to {\it stability}. It is also linked to {\it hyperbolicity}, a notion introduced by Smale~\cite{Sm}. More precisely, hyperbolic systems possess the shadowing property, which plays a crucial role in proving their stability. Pilyugin~\cite{P} demonstrated that structurally stable diffeomorphisms must satisfy a stronger form of the shadowing property. For a broader discussion on the significance of shadowing in both the qualitative theory of dynamical systems and numerical applications, we refer the reader to~\cite{Pa,Pi}.

Every continuous map on a compact metric space $X$ induces a continuous map $2^f$ (called the {\it induced map}) on the hyperspace $2^X$ of all nonempty closed subsets of $X$. If $X$ is connected, we consider the hyperspace $C(X)$ consisting of all nonempty closed and connected subsets of $X$. A naturally arising question is: what are the possible relations between the given (individual) dynamics on $X$ and the induced one (collective dynamics) on the hyperspace. Over the past few decades, various results have been obtained in this direction, yet this relationship remains largely unexplored and continues to be of significant interest. For instance,  it is known that certain dynamical properties of the system $(X,f)$ are preserved in the induced system $(2^X,2^f)$ (such as Li-Yorke chaos - see~\cite{GKLOPP} - and positive topological entropy - see~\cite{LR}). Conversely, some properties of $(2^X,2^f)$ also imply the same properties for $(X,f)$ (for example, transitivity - see~\cite{RF}). However, for some properties there is no implication in any direction (for example, neither Devaney chaos of $(X,f)$ implies Devaney chaos of $(2^X,2^f)$, nor Devaney chaos of $(2^X,2^f)$ implies Devaney chaos of $(X,f)$, see~\cite{GKLOPP}). Without attempting to be exhaustive, we mention just a few significant contributions in this area: Borsuk and Ulam~\cite{BU}, Bauer and Sigmund~\cite{BS}, Rom\'an-Flores~\cite{RF}, Banks~\cite{B}, Acosta, Illanes and M\'endez-Lango~\cite{AIM}.

It was proved in~\cite{FG} that $f$ has the shadowing property if and only if this is true for $2^f$. Additionally, if $C(f)$ has the shadowing property, the same is true for $f$~\cite{FG}. Morse-Smale diffeomorphisms are among the simplest dynamical systems, and they all possess the shadowing property. Regarding the shadowing property of $C(f)$, it was proved in~\cite{AB} that if $f:\mathbb{S}^1\to\mathbb{S}^1$ is a Morse-Smale diffeomorphism or if $f:\mathbb{S}^2\to\mathbb{S}^2$ is time-one-map of negative gradient system of Morse height function, then $C(f)$ does not satisfy the shadowing property. Additionally, recent results provide both positive and negative answers to this question in the context of transitive Anosov diffeomorphisms and dendrite monotone maps, see~\cite{CD}. More precisely, it is proved in~\cite{CD} that a wider class -- namely, continuum-wise hyperbolic homeomorphisms -- does not satisfy the shadowing property for the induced map on the hyperspace of continua, see also~\cite{ACS,ACCV,CCS,K1,K2}. However, the question whether $C(f)$ has the shadowing property when $f$ is Morse-Smale diffeomorphism remains an open question, even for $\mathbb{S}^n$.

Our contribution to this problem is the following negative result, which holds for any closed smooth manifold $M$.

\begin{thrm}\label{thrm:main}
For any time-one map $f$ of a negative gradient flow of a Morse function on a closed smooth manifold $M$ that satisfies the Morse-Smale condition, the induced homeomorphism $C(f)$ does not satisfy the shadowing property.\qed
\end{thrm}

\noindent{\bf Acknowledgements.} The authors thank the anonymous referee for many valuable comments and suggestions.

\section{Preliminaries}

Let us recall some notions and their properties that we will use in the proof.

\subsection{Shadowing}

Let $(X,d)$ be a compact metric space and $f:X\to X$ a continuous map. A {\it (positive) orbit} of a point $x$ is the set $\{f^n(x)\}_{n\in\mathbb{N}}$. If $f$ is {\it reversible}, i.e., $f$ is a homeomorphism, we can define a {\it full orbit} as the set $\{f^n(x)\}_{n\in\mathbb{Z}}$. We say that the set $A$ is {\it positively invariant} if $f(A)\subseteq A$.

\begin{defn} Let $\delta>0$. We say that the sequence $\{x_n\}_{n\in\mathbb{N}}$ (respectively $\{x_n\}_{n\in\mathbb{Z}}$) is $\delta$-{\it pseudo-orbit} if
$$d(x_{n+1},f(x_n))<\delta$$ for all $n\in\mathbb{N}$  (respectively $n\in\mathbb{Z}$). 
\end{defn}
One can also define a finite $\delta$-pseudo-orbit.

\begin{defn} Let $\varepsilon>0$. We say that a true orbit $\{f^n(x)\}_{n\in\mathbb{N}}$, of a point $x\in X$ $\varepsilon$-{\it shadows} a $\delta$-pseudo-orbit $\{x_n\}_{n\in\mathbb{N}}$, if
\begin{equation}\label{eq:shadowing}
d(f^n(x),x_n)<\varepsilon,
\end{equation}
for every $n\in\mathbb{N}$.
\end{defn}
If $f$ is reversible, we can define shadowing of a full $\delta$-pseudo-orbit $\{x_n\}_{n\in\mathbb{Z}}$, by requiring the condition~(\ref{eq:shadowing}) for every $n\in\mathbb{Z}$. 

In this paper we deal with Morse gradient system, which is reversible, so by shadowing we always assume the shadowing of a pseudo-orbit $\{x_n\}_{n\in\mathbb{Z}}$.

\begin{defn} We say that a reversible dynamical system $(X,f)$ has the {\it shadowing property} if for every $\varepsilon>0$ there exists $\delta>0$ such that for any $\delta$-pseudo-orbit $\{x_n\}_{n\in\mathbb{Z}}$ there exists a true orbit $\{f^n(x)\}_{n\in\mathbb{Z}}$ that $\varepsilon$-shadows it.
\end{defn}

\subsection{Hyperspaces and induced maps}

For a compact metric space $(X,d)$, the hyperspace $2^X$ is the set of all nonempty closed subsets of $X$.
The topology on $2^X$ is induced by the Hausdorff metric
$$d_H(A,B):=\inf\{\varepsilon>0\mid A\subset U_\varepsilon(B),\;B\subset U_\varepsilon(A)\},$$
where
\begin{equation}\label{eq:neighb}
U_\varepsilon(A):=\{x\in X\mid d(x,A)<\varepsilon\}.
\end{equation}
The obtained space $2^X$ is called a {\it hyperspace induced by $X$}, and it also turns out to be compact with respect to Hausdorff metric.

If $X$ is also connected (and so a continuum), then the set $C(X)$ of all connected and closed nonempty subsets of $X$ is also compact and connected. The set $C(X)$ is called the {\it hyperspace of subcontinua} of $X$.

If $f:X\to X$ is continuous, then it induces continuous maps
$$\begin{aligned}&2^f:2^X\to 2^X,\quad &2^f(A):=\{f(x)\mid x\in A\}\\
&C(f):C(X)\to C(X),\quad &C(f)(A):=\{f(x)\mid x\in A\}.
\end{aligned}$$
If $f$ is a homeomorphism, so are $2^X$ and $C(f)$.

In this paper we deal with the hyperspace $C(X)$.

We will use small latin letters $a$ for points in the initial space $X$, and capital latin letters $A$ for points in the induced hyperspace $C(X)$.

An open and a closed balls in $X$ will be denoted by $B(a,r)$ and $B[a,r]$. An open and a closed balls in $C(X)$ will be denoted by $B_H(A,r)$ and $B_H[A,r]$.

\subsection{Morse gradient systems}\label{subsec:Morse}

Fix a Riemannian metric $g$ on $M$. The {\it gradient vector field} of a smooth function $F:M\to\mathbb{R}$ induced by $g$ is the unique vector field $\nabla_gF$ satisfying
$$g_p((\nabla_gF)_p,\xi_p)=dF_p(\xi_p)$$ for every tangent vector $\xi_p\in T_pM$.

A critical point $p$ of a smooth function $F:M\to\mathbb{R}$ is said to be {\it non-degenerated} if the Hessian matrix at $p$ is non-degenerate (i.e., $p$ is a hyperbolic equilibrium of the gradient system induced by the vector field $\nabla F$; this definition does not depend on the choice of local coordinates or the Riemannian metric $g$).

Let $M$ be a smooth closed connected manifold and $F:M\to \mathbb{R}$ a smooth Morse function, meaning that all critical points of $F$ are non-degenerate. 
For a critical point $p$ of $F$, the {\it Morse index} of $p$ is defined as the dimension of the largest subspace of the tangent space $T_pM$ on which the Hessian is negative definite. It is equal to the number of negative eigenvalues of the Hessian matrix at $p$, which is independent of the choice of local coordinates. We denote by $m_F(p)$ the Morse index of $p$.

Let $\phi^t$ be the negative gradient flow defined by
$$
\frac{d\phi^t}{dt}(x)=-\nabla_g F(\phi^t(x)),\quad\phi^0=\mathrm{Id},
$$
 where the gradient $\nabla_g$ is induced by the metric $g$.
 
 For a critical point $p$ of $F$ define {\it unstable} and {\it stable manifold} of $p$ as:
$$W^u(p):=\{x\in M\mid \lim_{t\to-\infty}\phi^t(x)=p\},\quad W^s(p):=\{x\in M\mid \lim_{t\to+\infty}\phi^t(x)=p\}.$$
It is known that $W^u(p)$ and $W^s(p)$ are submanifolds of $M$ of dimension $m_F(p)$ and $\dim M-m_F(p)$ respectively (in fact they are diffeomorphic to $\mathbb{R}^{m_F(p)}$ and $\mathbb{R}^{\dim M-m_F(p)}$).
 
We say that that the pair $(F,g)$ satisfies {\it Morse-Smale condition} if for any two critical points $p$ and $q$, the manifolds $W^u(p)$ and $W^s(q)$ intersect transversally in $M$. Since the intersection of two submanifolds $N_1, N_2\subseteq M$ of codimensions $k_1$ and $k_2$ that intersect transversaly is either the empty set or a manifold of codimension $k_1+k_2$, we have that 
$$\mathcal{M}(p,q):=W^u(p)\cap W^s(q)$$ is either the empty set or a manifold of codimension 
$$\dim M-m_F(p)+(\dim M-(\dim M-m_F(q)))=\dim M-(m_F(p)-m_F(q)),$$ i.e., of the dimension $m_F(p)-m_F(q)$ (see also Subsection 2.2 in~\cite{AD}). The group $\mathbb{R}$ of translations in time acts on $\mathcal{M}(p,q)$ by:
$$\mathbb{R}\times\mathcal{M}(p,q)\ni(s,x) \mapsto\phi^s(x)$$ and this action is free for $p\neq q$. Therefore, the quotient $\mathcal{M}(p,q)/\mathbb{R}$ is a manifold, denoted by $\widehat{\mathcal{M}}(p,q)$ of dimension:
\begin{equation}\label{eq:dimM}
\dim\widehat{\mathcal{M}}(p,q)=m_F(p)-m_F(q)-1
\end{equation}
(see also Proposition 2.2.2 in~\cite{AD}). The space $\widehat{\mathcal{M}}(p,q)$ is called {\it the space of unparametrized trajectories}, where every trajectory $\gamma$ represents a single object, i.e., one point.

The time-one-map of Morse negative gradient equation satisfying Morse-Smale condition is a Morse-Smale diffeomorphism.

The manifold $\mathcal{M}(p,q)$ does not need to be closed, it can have a topological boundary that consists of ``broken trajectories", see~\cite{S}. In this paper we will use only one inclusion of this identification between the boundary of $\mathcal{M}(p,q)$ on the one hand, and the space of broken trajectories, on the other. To be precise, it is known that, for given pair of Morse trajectories $\alpha$ and $\beta$, satisfying Morse negative gradient equation
\begin{equation}\label{eq:Morse-grad-flow}
\frac{d\alpha}{dt}=-\nabla_gF(\alpha(t)),\quad\frac{d\beta}{dt}=-\nabla_gF(\beta(t))
\end{equation}
 and the boundary conditions:
$$\alpha(-\infty)=p,\;\alpha(+\infty)=\beta(-\infty)=q,\;\beta(+\infty)=r,$$
there exists a sequence $\gamma_n$ satisfying~(\ref{eq:Morse-grad-flow}) with boundary conditions
$$\gamma_n(-\infty)=p,\quad\gamma_n(+\infty)=r$$ that in some sense converges to the pair $(\alpha, \beta)$. The construction of this sequence is called {\it gluing}. In our proof the precise definition of this convergence is not relevant, we will use only the existence of this sequence. See Figure 1, or~\cite{S} for more details.

\vspace{5mm}

{\unitlength 0.8mm
\begin{picture}(150.00,50.00) \kern0pt\footnotesize
\put(30,0){

{\linethickness{1pt} \qbezier(30,50)(20,30)(0,20)
\qbezier(0,20)(22,12)(30,-10) \qbezier(11.5,16.5)(10.5,15.5)(13,13.3)\qbezier(13,13.3)(11,14)(10,12.5)
\qbezier(7,24.3)(8,24.5)(8,26.3)\qbezier(7,24.3)(9,24.5)(9,23)}

\qbezier(30,50)(29,35)(25,28)\qbezier(25,28)(20,21.5)(23,15)
\qbezier(23,15)(29,0)(30,-10)

\qbezier(30,50)(27,36)(20,30)\qbezier(20,30)(2,20)(17,15)
\qbezier(17,15)(25,12)(30,-10)

\put(14,21){\makebox(0,0)[l]{\normalsize{$\cdots$}}}

\put(25,21){\makebox(0,0)[l]{\normalsize{$\gamma_n$}}}
\qbezier(21,23)(22,23)(22,21)\qbezier(24,23)(22,23)(22,21)
\put(-11.7,-1){\qbezier(21,23)(22,23)(22,21)\qbezier(24,23)(22,23)(22,21)}

\put(11,35){\makebox(0,0)[l]{\normalsize{$\alpha$}}}
\put(11,7){\makebox(0,0)[l]{\normalsize{$\beta$}}}
\put(31,55){\makebox(0,0)[l]{\normalsize{$p$}}}
\put(33,-11){\makebox(0,0)[l]{\normalsize{$r$}}}
\put(-7,19){\makebox(0,0)[l]{\normalsize{$q$}}}
}

\end{picture}
\vspace{5mm}
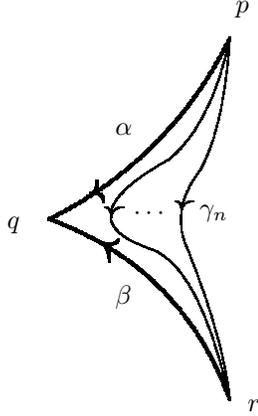
\begin{center}{\bf Figure 1:}\,{\rm Convergence to a broken trajectory}
\end{center}

\section{Proof of Theorem~\ref{thrm:main} }

In this section we prove our main result. For the reader's convenience, we will restate it.

\begin{thrm} Let $M$ be a smooth closed manifold and  $f=\phi^1$ the time-one map of negative Morse gradient flow~(\ref{eq:Morse-grad-flow}) which satisfies Morse--Smale condition. Then $C(f)$ does not have the shadowing property.
\end{thrm}

\noindent{\it Proof.} Suppose that $\dim M\ge 2$, since the case of $M=\mathbb{S}^1$ is done in~\cite{AB}. 

We need to find $\varepsilon>0$ and, for every $\delta$, a $\delta$-pseudo-orbit $\{X_n\}_{n\in\mathbb{Z}}$ that cannot be $\varepsilon$-shadowed. We will divide the proof in several steps. 

\vspace{3mm}

\noindent{\bf Step 1: construction of $X_0$ in pseudo-orbit. }

\vspace{2mm}

\begin{lem}
There exist two critical points $p$ and $q$ and two different solutions of negative gradient equation 
\begin{equation}\label{eq:2traj}
\gamma_i'(t)=-\nabla F(\gamma_i(t)),\quad i\in\{1,2\}\quad\mbox{ with}\quad \gamma_i(-\infty)=p, \,\gamma_i(+\infty)=q.
\end{equation}
\end{lem}

\noindent{\it Proof.} We distinguish between two cases.

Assume that there exists a critical point $r$ of Morse index $0<m_F(r)<\dim M$. Let $\alpha$ and $\beta$ be any two gradient curves with $\alpha(+\infty)=r$ and $\beta(-\infty)=r$. We can prove that these two curves exist by using for example the Hartman-Grobman theorem. Indeed, since $F$ is Morse, $r$ is hyperbolic critical point of $\nabla F$ (meaning that $L:=-d\left(\nabla F(r)\right)$ is a hyperbolic matrix). The Hartman-Grobman theorem says that locally, the dynamical system induced by the differential equation~(\ref{eq:Morse-grad-flow}) is equivalent to the dynamical system defined by linearized system:
$$ \frac{d\psi^t}{dt}(x)=L\cdot\psi^t(x),\quad\psi^0=\mathrm{Id}.$$ Since $0<m_F(r)<\dim M$, the symmetric matrix $L$ has both positive and negative eigenvalues, which implies that there exist at least one non-constant trajectory, $\beta$, which originates at $r$, and at least one non-constant trajectory, $\alpha$, that ends at $r$.

Let $p:=\alpha(-\infty)$ and $q:=\beta(+\infty)$. Since there exists a non-constant trajectory from $p$ to $r$, the $\mathbb{R}$ action on $\mathcal{M}(p,r)$ is non-trivial. Consequently, the set $\widehat{\mathcal{M}}(p,r)$ is non-empty, and thus has dimension at least zero. From~(\ref{eq:dimM}) we compute:
$$m_F(p)-m_F(r)-1\ge 0\quad\Rightarrow\quad m_F(p)\ge m_F(r)+1$$ and similarly $m_F(r)\ge m_F(q)+1$. Therefore $m_F(p)\ge m_F(q)+2$, so the dimension of the manifold $\mathcal{M}(p,q)=W^u(p)\cap W^s(q)$ is at least two, if it is nonempty. From the discussion in Subsection~\ref{subsec:Morse}, we can conclude that $W^u(p)$ and $W^s(q)$ must intersect since there exists a broken trajectory $(\alpha,\beta)$. Therefore there exist infinitely many trajectories satisfying~(\ref{eq:2traj}).

In the second case there are no critical points of Morse index $0<m_F(r)<\dim M$, we can take any gradient trajectory to be $\gamma_1$ and define $p:=\gamma_1(-\infty)$ and $q:=\gamma_1(+\infty)$. Then the dimension of $\mathcal{M}(p,q)$ is equal to $m_F(p)-m_F(q)=\dim M\ge 2$, so the cardinality of $\widehat{\mathcal{M}}(p,q)$ is infinity, if it is nonempty. Since $\gamma_1$ belongs to $\widehat{\mathcal{M}}(p,q)$, we find $\gamma_2$ with the same boundary condition. \qed

\vspace{3mm}

Now fix $\gamma_1$ and $\gamma_2$ satisfying~(\ref{eq:2traj}) and define 
$$X_0:=\gamma_1\cup\gamma_2\in C(M).$$

\vspace{3mm}

\noindent{\bf Step 2: construction of $\varepsilon$. }

\vspace{2mm}
Denote by $b:=F(p)$ and $a:=F(q)$. Since $F$ decreases along its negative gradient flow, we have $a<b$. Choose $a_1,b_1\in\mathbb{R}$ such that $a<a_1<b_1<b$. Since the sets
$$A_1:=\{\gamma_1(t)\mid a_1\le F(\gamma_1(t))\le b_1\}\quad\mbox{and}\quad A_2:=\{\gamma_2(t)\mid a_1\le F(\gamma_2(t))\le b_1\}$$ are compact and disjoint, there exists $\varepsilon>0$ such that
$$\overline{U_\varepsilon(A_1)}\cap\overline{ U_\varepsilon(A_2)}=\emptyset,$$ where $U_\varepsilon(\cdot)$ is defined in~(\ref{eq:neighb}), see Figure 2.

We can decrease $a_1\in(a,b_1)$ (note that this may result in decreasing $\varepsilon$) if necessary, to obtain 
\begin{equation}\label{eq:no-jump}
x\in B_H[X_0,\varepsilon]\cap\mathcal{M}(p,q)\cap\{F\ge b_1\}\quad\Rightarrow\quad F(f(x))> a_1.
\end{equation}
 This is possible to do since for every such $x$ it holds $F(f(x))>a$ and the set $B_H[X_0,\varepsilon]\cap\mathcal{M}(p,q)\cap\{F\ge b_1\}$ is compact in $\mathcal{M}(p,q)$.

We will also decrease $\varepsilon$ if necessary to get the following implication:
\begin{equation}\label{eq:stays}
x\in U_\varepsilon(A_1)\;\Rightarrow\; f(x)\notin U_\varepsilon(A_2),\quad\mbox{and}\quad
x\in U_\varepsilon(A_2)\;\Rightarrow\; f(x)\notin U_\varepsilon(A_1)
\end{equation} 
(this can be done since it holds for every $x\in A_1$ or $x\in A_2$, and $A_1$ and $A_2$ are compact).

\vspace{1cm}

{\unitlength 0.7mm
\begin{picture}(150.00,50.00) \kern0pt\footnotesize
\put(0,0){

{\linethickness{1pt} 
\qbezier(30,50)(0,20)(30,-10)\qbezier(30,50)(60,20)(30,-10)}

\qbezier(-5,45)(32,39)(65,45)
\qbezier(-5,-5)(32,-1)(65,-5)
\begin{color}{red} {\linethickness{0.8pt}\qbezier[70](40,42)(57,20)(40,-3)
\qbezier[70](33,42)(50,20)(33,-3)}
\put(45,37){\makebox(0,0)[l]{\normalsize{$U_\varepsilon(A_2)$}}}

\end{color}
}

\begin{color}{green} {\linethickness{0.8pt}\qbezier[70](27,42)(10,20)(27,-3)
\qbezier[70](20,42)(3,20)(20,-3)
}
\put(0,37){\makebox(0,0)[l]{\normalsize{$U_\varepsilon(A_1)$}}}
\end{color}

\put(81,0){

{\linethickness{1pt} 

\qbezier(30,50)(0,20)(30,-10)\qbezier(30,50)(60,20)(30,-10)
}
\qbezier(-5,45)(32,39)(65,45)\qbezier(-5,-5)(32,-1)(65,-5)

\begin{color}{blue}
\put(0.3,0){\qbezier(32,42)(42,52)(42,52)
\put(4,0){\qbezier(32,42)(42,52)(42,52)}
\put(8,0){\qbezier(32,42)(42,52)(42,52)}
\put(0.6,0){\put(12,0.7){\qbezier(32,42)(42,52)(42,52)}
\put(16,1){\qbezier(32,42)(42,52)(42,52)}}}

\put(-18.4,0){\qbezier(32,43)(42,53)(42,53)
\put(4,0){\qbezier(31.5,42.5)(41.5,52.5)(41.5,52.5)}
\put(8,0){\qbezier(31.5,42)(42,52.5)(41.5,52.5)}
\put(0.6,0){\put(12,0.7){\qbezier(31.5,41.5)(42,52)(42,52)}
\put(16,1){\qbezier(31,41)(41,51)(41,51)}}}

\put(-34.7,2){\qbezier(32,43)(42,53)(42,53)
\put(4,0){\qbezier(31.5,42.5)(41.5,52.5)(41.5,52.5)}
\put(8,0){\qbezier(31.5,42)(42,52.5)(41.5,52.5)}
\put(0.6,0){\put(12,0.1){\qbezier(31.5,41.5)(42,52)(42,52)}}}
\put(56,47){\makebox(0,0)[l]{\normalsize{$A$}}}
\end{color}

\begin{color}{brown}
\put(0,-55.8){\qbezier(32,43)(42,53)(42,53)
\put(4,0){\qbezier(31.5,42.5)(41.5,52.5)(41.5,52.5)}
\put(8,0){\qbezier(31.5,42)(42,52.5)(41.5,52.5)}
\put(-0.5,0){\put(12,0.1){\qbezier(31.5,41.5)(42,52)(41.7,52)}}
\put(56,47){\makebox(0,0)[l]{\normalsize{$B$}}}
}
\put(-16.2,-55){\put(0,-0.5){\qbezier(32,43)(42,53)(42,53)}
\put(4,-0.5){\qbezier(31.5,42.5)(41.5,52.5)(41.5,52.5)}
\put(8,-0.5){\qbezier(31.5,42)(42,52.5)(41.5,52.5)}
\put(-0.5,-0.5){\put(12,0.1){\qbezier(31.5,41.5)(42,52)(41.7,52)}}}

\put(-32,-55){\put(-0.5,-1.5){\qbezier(32,43)(42,53)(42,53)}
\put(4,-0.5){\qbezier(31.5,42.5)(41.5,52.5)(41.5,52.5)}
\put(8,-0.5){\qbezier(31.5,42)(42,52.5)(41.5,52.5)}
\put(-0.5,-0.5){\put(12,0.1){\qbezier(31.5,41.5)(42,52)(41.7,52)}}}
\end{color}
}

\put(87,30){\makebox(0,0)[l]{\normalsize{$X_0$}}} 

\end{picture}}
\vspace{13mm}
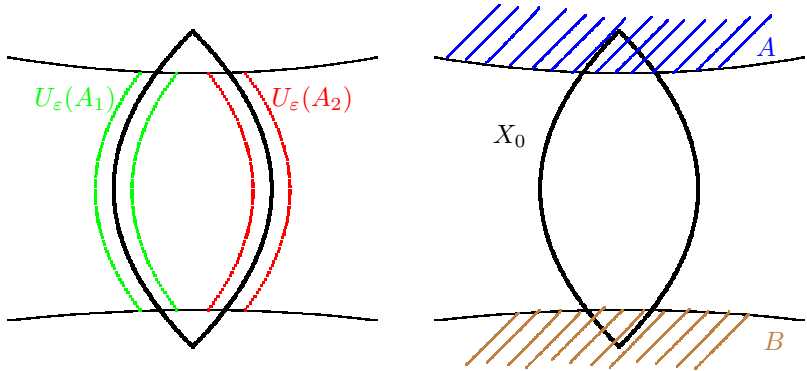
\begin{center}{\bf Figure 2:}\,{\rm On the left: sets $A_1$ and $A_2$ and their $\varepsilon$-neighborhood. On the right: sets $A$ and $B$}
\end{center}
\vspace{8mm}

\noindent{\bf Step 3: construction of $\delta$-pseudoorbit $\{X_n\}_{n\in\mathbb{Z}}$. }
\vspace{2mm}

Now we use the idea from the proof of Theorem A in~\cite{AB}. Recall that $X_0:=\gamma_1\cup\gamma_2\in C(M)$. For given $\delta>0$, choose $M>0$ such that
$$d_H(\gamma_1((-\infty,M])\cup\gamma_2((-\infty,M]),X_0)<\delta, \quad d_H(\gamma_1([-M,\infty))\cup\gamma_2([-M,\infty)),X_0)<\delta.$$ Define
\begin{itemize}
\item $X_{1}:=\gamma_1([-M,\infty))\cup\gamma_2([-M,\infty))$
\item $X_{-1}:=\gamma_1((-\infty,M])\cup\gamma_2((-\infty,M])$
\item $X_i:=C(f)^i(X_1)$, for $i>1$
\item $X_i:=C(f)^i(X_{-1})$, for $i<-1$, 
\end{itemize}
see Figure 3.

We have constructed a $\delta$-pseudo-orbit $\{X_n\}_{n\in\mathbb{Z}}\subset C(f)$. Note that
\begin{equation}\label{eq:lim-of-x_n}
X_n\to 
\begin{cases} \{q\}, & n\to \infty \\
\{p\}, & n\to -\infty.
\end{cases}
\end{equation}

\vspace{1cm}

{\unitlength 0.8mm
\begin{picture}(150.00,50.00) \kern0pt\footnotesize

\put(-20,0){
{\linethickness{1pt} \qbezier(30,50)(0,20)(30,-10)\qbezier(30,50)(60,20)(30,-10)}
\put(33,30){\makebox(0,0)[l]{\normalsize{$X_0$}}} 
}

\put(35,0){
{\linethickness{0.5pt} \qbezier(30,50)(0,20)(30,-10)\qbezier(30,50)(60,20)(30,-10)}
\begin{color}{red}
\linethickness{1.1pt}{ \qbezier(25,45)(1,20)(29,-10)
\qbezier(33,45)(57,20)(29,-10)}
\put(33,30){\makebox(0,0)[l]{\normalsize{$X_{-1}$}}} 

\end{color}}
\put(90,0){
{\linethickness{0.5pt} \qbezier(30,50)(0,20)(30,-10)\qbezier(30,50)(60,20)(30,-10)}
\begin{color}{green}
{\linethickness{1.1pt} \qbezier(29,50)(1,20)(25,-5)
\qbezier(29,50)(57,20)(33,-5)}

\put(33,30){\makebox(0,0)[l]{\normalsize{$X_1$}}} 
\end{color}
}

\end{picture}}
\vspace{5mm}
\begin{center}{\bf Figure 3:}\,{\rm $\delta$-pseudo-orbit $X_n$}
\end{center}

\vspace{1cm}

\noindent{\bf Step 4: the end of the proof. }

\vspace{2mm}

Suppose that there exists $K\in C(M)$ that $\varepsilon$-shadows $\{X_n\}_{n\in\mathbb{Z}}$. Denote by
$$A:=U_\varepsilon(X_0)\cap \{x\mid F(x)\ge b_1\}\quad B:=U_\varepsilon(X_0)\cap \{x\mid F(x)\le a_1\},$$
see Figure 2.
We conclude from~(\ref{eq:lim-of-x_n}) that there exists $n_0\in\mathbb{N}$ such that
$f^{-n_0}(K)\subset A$. Denote by $K_0:=f^{-n_0}(K)$. Since $K$ is connected, so it is $K_0$.

For any point $x\in K_0$ there exists $n\ge 1$ such that $f^n(x)\in U_\varepsilon(A_1)\cup U_\varepsilon(A_2)$. Indeed, choose a minimal $k\ge 1$ such that $f^k(x)\notin A$. It follows from~(\ref{eq:no-jump}) that $f^k(x)\notin B$, therefore, since $K$ $\varepsilon$-shadows $\{X_n\}_{n\in\mathbb{Z}}$, $f^k(x)$ must be contained in $U_\varepsilon(A_1)\cup U_\varepsilon(A_2)$. Denote by $k_x$ this minimal number $k$, depending on $x$, such that $f^{k_x}(x)\in U_\varepsilon(A_1)\cup U_\varepsilon(A_2)$.

Define the following function $\varphi:K_0\to\{0,1\}$:
$$\varphi(x)=\begin{cases}0,&f^{k_x}(x) \in U_\varepsilon(A_1),\\
1,& f^{k_x}(x)\in U_\varepsilon(A_2).\end{cases}$$
This function is continuous. To see this, suppose that $x_0\in K_0$ and $\varphi(x_0)=0$ (the other case is treated in the same way). Let $k$ be the smallest integer such that $f^k(x_0)\in U_\varepsilon(A_1)$.

Since $f^k$ is continuous, there exists a neighbourhood $U_{x_0}$ of $x_0$ such that 
$$f^k(U_{x_0}\cap K_0)\subset U_\varepsilon(A_1).$$ It follows from~(\ref{eq:stays}) that for every $y\in U_{x_0}$ and $j\in\{0,\ldots,k-1\}$ we have $f^j(y)\notin U_\varepsilon(A_2)$. This implies that $f^{k_y}(y)\in U_\varepsilon(A_1)$, i.e., $\varphi(y)=0$ for all $y\in U_{x_0}$. Therefore, $\varphi$ is continuous at $x_0$.

Since $K_0$ is connected, we conclude that $\varphi$ is constant, suppose $\varphi=0$. This means that every point $x\in K_0$ enters $U_\varepsilon(A_1)$. From~(\ref{eq:stays}) it follows that if $f^k(x)\in U_\varepsilon(A_1)$
must imply either $f^{k+1}(x)\in U_\varepsilon(A_1)$ or $f^{k+1}(x)\in B$. Since $F$ decreases along the orbits of $f$, the set $B$ is $f$-positive invariant. This means that every point $x$ that enters $B$, cannot enter $U_\varepsilon(A_2)$, implying $f^n(K_0)\cap U_\varepsilon(A_2)=\emptyset$, for every $n$, so $K$ does not $\varepsilon$-shadow $\{X_n\}_{n\in\mathbb{Z}}$.\qed

 \end{document}